\documentclass[12pt,a4paper]{article}

\usepackage[T1]{fontenc}
\usepackage{lmodern}
\usepackage{amsmath,amssymb,amsthm}
\usepackage[margin=28mm]{geometry}
\usepackage{microtype}
\usepackage{needspace}
\usepackage{cite}
\usepackage[hidelinks]{hyperref}

\hypersetup{
  pdftitle={Local superderivations of Block Lie superalgebras},
  pdfauthor={Jakhongir Ruzimboev and Bakhtiyor B. Yusupov},
  pdfkeywords={Block Lie superalgebra, local superderivation, homogeneous superderivation}
}

\newtheorem{theorem}{Theorem}[section]
\newtheorem{proposition}[theorem]{Proposition}
\newtheorem{lemma}[theorem]{Lemma}
\theoremstyle{definition}
\newtheorem{definition}[theorem]{Definition}

\numberwithin{equation}{section}
\title{\bfseries Local superderivations of Block Lie superalgebras $S(q)$}
\author{
Jakhongir Ruzimboev\thanks{New Uzbekistan University, Tashkent, Uzbekistan;
Tashkent International University of Financial Management
and Technologies, Tashkent, Uzbekistan.
Email: \href{mailto:jahongirrozimboyev5@gmail.com}{\nolinkurl{jahongirrozimboyev5@gmail.com}}.}
\and
Bakhtiyor Yusupov\thanks{V.I. Romanovskiy Institute of Mathematics,
Uzbekistan Academy of Sciences, Tashkent, Uzbekistan;  Department of Applied Mathematics and Mathematical Physics, Urgench State University, Urgench, Uzbekistan.
Email: \href{mailto:baxtiyor_yusupov_93@mail.ru}{\nolinkurl{baxtiyor_yusupov_93@mail.ru}}.}
}
\date{}

\begin{document}

\maketitle

\begin{abstract}
We prove that every local superderivation of the
Block Lie superalgebra $S(q)$ is a superderivation for every
$q\in\mathbb C$, with both basis indices ranging over $\mathbb Z$.
The superderivationing derivations have the same parity as the local map.
The proof uses the known derivation formulas and coefficient
comparisons.
\end{abstract}



\section{Introduction}

Local derivations were introduced by Kadison~\cite{Kadison}
and independently studied by Larson and Sourour~\cite{LarsonSourour}.
A local derivation is a linear map that agrees at each element
with a derivation which may depend on that element.
The main question is whether this condition forces the map
to be a derivation on the whole algebra.
A related notion is that of a $2$-local derivation, introduced
by \v{S}emrl~\cite{Semrl}.
For every pair of elements, a single derivation must agree
with the map at both elements.
A $2$-local derivation is not assumed to be linear.

For finite-dimensional Lie algebras, Ayupov and
Kudaybergenov~\cite{AyupovKudaybergenov} proved that every
local derivation of a semisimple Lie algebra is a derivation.
They also gave examples of nilpotent Lie algebras for which
this conclusion fails.
For infinite-dimensional Lie algebras, Chen, Zhao and
Zhao~\cite{ChenZhaoZhao} proved the corresponding result
for Witt algebras.
Local derivations of Schr\"odinger algebras were studied
by Alauadinov and Yusupov~\cite{AlauadinovYusupovSchrodinger}.

For Lie superalgebras, Chen, Wang and Nan~\cite{ChenWangNan}
proved that every local superderivation of a basic classical
Lie superalgebra, except $A(1,1)$, is a superderivation.
Yuan, Chen and Cao~\cite{YuanChenCao} obtained the corresponding
result for Cartan type Lie superalgebras.
Local superderivations of super Schr\"odinger algebras and the
modified Witt-type Lie superalgebra were studied in
\cite{AlauadinovYusupovSuperSchrodinger,RuzimboevBoranbaev},
where every local superderivation was shown to be a
superderivation.

Block introduced a class of infinite-dimensional Lie algebras
in~\cite{Block}.
Biderivations and transposed Poisson structures on Block Lie
algebras and related superalgebras were studied in
\cite{LiuGuoZhao,KaygorodovKhrypchenko}.
Wu and Yang~\cite{WuYang} proved that every local derivation
of the Block-type Lie algebra with $q=1$ and nonnegative second
basis index is a derivation.
Chen and Guo~\cite{ChenGuoBlockTwoLocal} proved the corresponding
$2$-local result for a related Block-type Lie algebra with
nonnegative second basis index.
For the Virasoro-like algebra, the derived algebra of
$\mathcal B(0)$, Tang, Xiao and Wang
\cite{TangXiaoWangVirasoroLike} proved that every $2$-local
derivation is a derivation.
For $\mathcal B(q)$ with both basis indices in $\mathbb Z$,
every local or $2$-local derivation was shown to be a
derivation for all $q\in\mathbb C$
in~\cite{RuzimboevYusupovBlock}.

Xia~\cite{XiaBlockSuper} studied the structure of the
Ramond--Block and Neveu--Schwarz--Block Lie superalgebras.
For positive integral parameters and nonnegative second
basis indices, Su, Yue and Zhu~\cite{SuYueZhu} proved that
every $2$-local superderivation of these algebras is a
superderivation.
The homogeneous derivations and quasi-derivations of the
Block Lie superalgebras $S(q)$, with both basis indices in
$\mathbb Z$ and arbitrary complex parameter $q$, were
described in~\cite{KaygorodovKhudoyberdiyevShermatovaSuper}.
The corresponding Block Lie algebras $\mathcal B(q)$
were studied in~\cite{KaygorodovKhudoyberdiyevShermatova}.

In this paper, we prove that every even or odd local
superderivation of $S(q)$ is a superderivation for all
$q\in\mathbb C$.
The superderivations used in the local condition are required
to have the same parity as the map.
Locality is imposed on all elements, including sums of even
and odd elements.
We treat the even and odd cases separately and distinguish
between $q\ne0$ and $q=0$.
The proofs use the known description of homogeneous
superderivations, explicit bracket computations, and
coefficient comparisons on basis elements and their finite sums.
\section{Preliminaries}\label{sec:preliminaries}

All vector spaces and algebras are over $\mathbb C$.
All sums indexed by basis elements are finite.
We write $\delta_{a,b}$ for the Kronecker symbol and
$\mathbb Z_2=\{\bar0,\bar1\}$ for the group of parities.

Let $S=S_{\bar0}\oplus S_{\bar1}$ be a Lie superalgebra.
For a homogeneous element $x$, its parity is denoted by $|x|$.
A linear map $D$ has parity $\theta\in\mathbb Z_2$ if
$D(S_{\alpha})\subseteq S_{\alpha+\theta}$ for every
$\alpha\in\mathbb Z_2$.

\begin{definition}
A linear map $D:S\to S$ of parity $\theta$ is a
\emph{superderivation} if
\begin{equation*}
D([x,y])=[D(x),y]+(-1)^{\theta|x|}[x,D(y)]
\end{equation*}
for all homogeneous $x,y\in S$.
The space of these maps is denoted by
$\operatorname{Der}_{\theta}(S)$.
\end{definition}

For $u\in S_{\theta}$, the map
$\operatorname{ad}u$, defined by
$(\operatorname{ad}u)(x)=[u,x]$, is an inner superderivation
of parity $\theta$.

\begin{definition}\label{def:local-super}
A linear map $\Delta:S\to S$ of parity $\theta$ is called a
\emph{local superderivation of parity $\theta$} if, for every
$x\in S$, there exists
$D_x\in\operatorname{Der}_{\theta}(S)$ such that
$\Delta(x)=D_x(x)$.
We denote the space of these maps by
$\operatorname{LDer}_{\theta}(S)$.
The map is called even if $\theta=\bar0$ and odd if
$\theta=\bar1$.
\end{definition}

In Definition~\ref{def:local-super}, locality is required
at every element, including sums of even and odd elements.
If $x=x_{\bar0}+x_{\bar1}$, with
$x_{\bar0}\in S_{\bar0}$ and $x_{\bar1}\in S_{\bar1}$,
then linearity and comparison of parity components give
\begin{equation*}
D_x(x_{\bar0})=\Delta(x_{\bar0}),
\qquad
D_x(x_{\bar1})=\Delta(x_{\bar1}).
\end{equation*}
Here the same derivation $D_x$ works at both components.
Also, if $\Delta\in\operatorname{LDer}_{\theta}(S)$ and
$D_0\in\operatorname{Der}_{\theta}(S)$, then
$\Delta-D_0\in\operatorname{LDer}_{\theta}(S)$:
at $x$, use the derivation $D_x-D_0$.

For $q\in\mathbb C$, the Block Lie superalgebra
$S(q)=S(q)_{\bar0}\oplus S(q)_{\bar1}$ has even basis
$\{L_{m,i}:m,i\in\mathbb Z\}$ and odd basis
$\{G_{m,i}:m,i\in\mathbb Z\}$, with brackets
\begin{equation}\label{eq:super-brackets}
\begin{aligned}
{}[L_{m,i},L_{n,j}]
&=\bigl(n(i+q)-m(j+q)\bigr)L_{m+n,i+j},\\
[L_{m,i},G_{n,j}]
&=\bigl(n(i+q)-m(j+q/2)\bigr)G_{m+n,i+j},\\
[G_{m,i},G_{n,j}]
&=2qL_{m+n,i+j}.
\end{aligned}
\end{equation}
All indices range over $\mathbb Z$.
The remaining brackets follow from
$[x,y]=-(-1)^{|x||y|}[y,x]$.
In particular, the bracket is symmetric on the odd part.
The even subalgebra is the Block Lie algebra $\mathcal B(q)$.
For integral $q$, the element $L_{0,-q}$ is central.
When $q=0$, the element $G_{0,0}$ is central as well.

We use the notation for the following superderivations from
\cite{KaygorodovKhudoyberdiyevShermatovaSuper}.
For $X\in\{L,G\}$, put
$D_1(X_{m,i})=iX_{m,i}$ and
$D_2(X_{m,i})=mX_{m,i}$.
Both maps are even.
When $q\ne0$, the map $D_2=q^{-1}\operatorname{ad}L_{0,0}$
is inner.

For $q\in\mathbb Z\setminus\{0\}$, define the even maps
\begin{equation*}
\begin{aligned}
D_3(L_{m,i})&=mL_{m,i-q},
&
D_3(G_{m,i})&=mG_{m,i-q},\\
D_4(L_{m,i})&=(i+q)L_{m,i-q},
&
D_4(G_{m,i})&=(i+q/2)G_{m,i-q}.
\end{aligned}
\end{equation*}
For $q=0$, let $E_L(L_{0,0})=L_{0,0}$ and
$E_G(G_{0,0})=G_{0,0}$, with each map zero on all other
basis vectors.
We also use the even map $T$ given by
$T(L_{m,i})=0$ and $T(G_{m,i})=G_{m,i}$.

For $q\in2\mathbb Z\setminus\{0\}$, define the odd map
$D_6$ by $D_6(G_{0,-3q/2})=L_{0,-q}$ and set its value
to zero on every other basis vector.
For $q=0$, define the odd maps
\begin{equation*}
\begin{aligned}
O_1(L_{m,i})&=iG_{m,i},
&
O_1(G_{m,i})&=0,\\
O_2(L_{m,i})&=mG_{m,i},
&
O_2(G_{m,i})&=0,\\
\Theta(L_{m,i})&=0,
&
\Theta(G_{m,i})&=L_{m,i}.
\end{aligned}
\end{equation*}
Finally, let $O_5(L_{0,0})=G_{0,0}$ and
$\Theta_0(G_{0,0})=L_{0,0}$, with each map zero on all
other basis vectors.

The description of homogeneous superderivations below
is a known result.

\begin{proposition}[
\cite{KaygorodovKhudoyberdiyevShermatovaSuper}, Theorem~8]
\label{prop:super-derivations}
For every $q\in\mathbb C$,
\begin{equation}\label{eq:even-derivation-form}
\operatorname{Der}_{\bar0}(S(q))
=
\begin{cases}
\operatorname{ad}(S(q)_{\bar0})\oplus\mathbb C D_1,
&q\notin\mathbb Z,\\
\operatorname{ad}(S(q)_{\bar0})
 \oplus\operatorname{span}_{\mathbb C}\{D_1,D_3,D_4\},
&q\in\mathbb Z\setminus\{0\},\\
\operatorname{ad}(S(0)_{\bar0})
 \oplus\operatorname{span}_{\mathbb C}\{D_1,D_2,E_L,E_G,T\},
&q=0.
\end{cases}
\end{equation}
For odd superderivations,
\begin{equation}\label{eq:odd-derivation-form}
\operatorname{Der}_{\bar1}(S(q))
=
\begin{cases}
\operatorname{ad}(S(q)_{\bar1}),
&q\notin2\mathbb Z,\\
\operatorname{ad}(S(q)_{\bar1})\oplus\mathbb C D_6,
&q\in2\mathbb Z\setminus\{0\},\\
\operatorname{ad}(S(0)_{\bar1})
 \oplus\operatorname{span}_{\mathbb C}\{O_1,O_2,O_5,\Theta,\Theta_0\},
&q=0.
\end{cases}
\end{equation}
All sums are direct sums of vector spaces.
\end{proposition}

For a local superderivation $\Delta$ of fixed parity, each
superderivation $D_x$ has the corresponding form in
\eqref{eq:even-derivation-form} or
\eqref{eq:odd-derivation-form}.
Its inner element and scalar coefficients may depend on $x$.
For each fixed $x$, the inner element is a finite linear
combination of basis vectors.

\section{Local even superderivations}
\label{sec:even-local-super}

We study local even superderivations of $S(q)$, treating
the three parameter cases separately.
The following basis formula will be used in both
nonzero parameter cases.

\begin{lemma}\label{lem:even-basis-action}
Let $q\ne0$ and
$\Delta\in\operatorname{LDer}_{\bar0}(S(q))$.
There are finite families of scalars
$p_{r,k},s_{r,k},c_{r,k}\in\mathbb C$ such that
\begin{equation}\label{eq:even-basis-action}
\begin{aligned}
\Delta(L_{m,i})
&=\sum_{r,k}(mp_{r,k}+is_{r,k}+2c_{r,k})L_{m+r,i+k},\\
\Delta(G_{m,i})
&=\sum_{r,k}(mp_{r,k}+is_{r,k}+c_{r,k})G_{m+r,i+k}
\end{aligned}
\end{equation}
for all $m,i\in\mathbb Z$.
\end{lemma}

\begin{proof}
Write
\begin{equation*}
\Delta(L_{m,i})
=\sum_{r,k}A_{r,k}(m,i)L_{m+r,i+k},
\qquad
\Delta(G_{m,i})
=\sum_{r,k}B_{r,k}(m,i)G_{m+r,i+k}.
\end{equation*}
For an odd element $g$, apply locality to $x=[g,g]+g$.
Comparison of parity components gives
$D_x(g)=\Delta(g)$ and $D_x([g,g])=\Delta([g,g])$.
Since $D_x$ is even and the bracket is symmetric on the odd
part, we obtain
\begin{equation*}
\Delta([g,g])
=[D_x(g),g]+[g,D_x(g)]
=2[\Delta(g),g].
\end{equation*}
Applying this identity to $g+h$, where $g,h$ are odd,
and using linearity to subtract the identities for $g$ and $h$,
gives
\begin{equation*}
\begin{aligned}
2\Delta([g,h])
&=\Delta([g+h,g+h])-\Delta([g,g])-\Delta([h,h])\\
&=2[\Delta(g),h]+2[g,\Delta(h)].
\end{aligned}
\end{equation*}
Take $g=G_{m,i}$ and $h=G_{n,j}$.
Using \eqref{eq:super-brackets} and $q\ne0$, comparison
of the coefficient of $L_{m+n+r,i+j+k}$ gives
\begin{equation}\label{eq:even-coefficient-relation}
A_{r,k}(m+n,i+j)=B_{r,k}(m,i)+B_{r,k}(n,j).
\end{equation}

Fix $r,k$, put $c_{r,k}=B_{r,k}(0,0)$, and define
$F_{r,k}(m,i)=B_{r,k}(m,i)-c_{r,k}$.
Setting $(n,j)=(0,0)$ in
\eqref{eq:even-coefficient-relation} gives
$A_{r,k}(m,i)=B_{r,k}(m,i)+c_{r,k}$.
Substituting this expression back into
\eqref{eq:even-coefficient-relation}, we obtain
\begin{equation*}
F_{r,k}(m+n,i+j)=F_{r,k}(m,i)+F_{r,k}(n,j).
\end{equation*}
Thus $F_{r,k}$ is additive on $\mathbb Z^2$.
Since $(m,i)=m(1,0)+i(0,1)$, it follows that
$F_{r,k}(m,i)=mp_{r,k}+is_{r,k}$, where
$p_{r,k}=F_{r,k}(1,0)=B_{r,k}(1,0)-c_{r,k}$ and
$s_{r,k}=F_{r,k}(0,1)=B_{r,k}(0,1)-c_{r,k}$.
Consequently,
\begin{equation*}
\begin{aligned}
A_{r,k}(m,i)&=mp_{r,k}+is_{r,k}+2c_{r,k},\\
B_{r,k}(m,i)&=mp_{r,k}+is_{r,k}+c_{r,k}.
\end{aligned}
\end{equation*}
The three families have finite support, since they are
obtained from the finite coefficient families of
$\Delta(G_{0,0})$, $\Delta(G_{1,0})$ and $\Delta(G_{0,1})$.
This proves \eqref{eq:even-basis-action}.
\end{proof}

\subsection{\texorpdfstring{The case $q\notin\mathbb Z$}
{The case q not in Z}}
\label{subsec:even-nonintegral}

Assume $q\notin\mathbb Z$, so $q\ne0$.
By Proposition~\ref{prop:super-derivations}, each even
superderivation has the form $D_x=\operatorname{ad}u_x+\alpha_xD_1$,
where $u_x\in S(q)_{\bar0}$ and $\alpha_x\in\mathbb C$.

\begin{lemma}\label{lem:even-local-nonintegral}
If $q\notin\mathbb Z$, every local even superderivation
of $S(q)$ is an even superderivation.
\end{lemma}

\begin{proof}
Let $\Delta\in\operatorname{LDer}_{\bar0}(S(q))$.
Choose an even superderivation $D_0$ at $L_{0,0}+G_{0,1}$ and put
$F=\Delta-D_0$.
Then $F$ is local even and $F(L_{0,0})=F(G_{0,1})=0$.
Apply Lemma~\ref{lem:even-basis-action} to $F$.
The equality $F(L_{0,0})=0$ forces $c_{r,k}=0$,
and $F(G_{0,1})=0$ then forces $s_{r,k}=0$.
Consequently,
\begin{equation}\label{eq:even-nonintegral-normalized}
\begin{aligned}
F(L_{m,i})&=m\sum_{r,k}p_{r,k}L_{m+r,i+k},\\
F(G_{m,i})&=m\sum_{r,k}p_{r,k}G_{m+r,i+k}.
\end{aligned}
\end{equation}

Apply locality for $F$ at $y=L_{0,0}+G_{1,0}$.
Write its superderivation as
$D_y=\operatorname{ad}u+\alpha D_1$, with
$u=\sum_{r,k}a_{r,k}L_{r,k}$.
Comparison of parity components gives
$D_y(L_{0,0})=0$ and $D_y(G_{1,0})=F(G_{1,0})$.
Since $D_1(L_{0,0})=0$, we get
\begin{equation*}
0=D_y(L_{0,0})
=\sum_{r,k}a_{r,k}[L_{r,k},L_{0,0}]
=-q\sum_{r,k}r a_{r,k}L_{r,k}.
\end{equation*}
Thus $a_{r,k}=0$ for $r\ne0$, and
$u=\sum_k a_{0,k}L_{0,k}$.
Also, $D_1(G_{1,0})=0$, so
\begin{equation*}
F(G_{1,0})=D_y(G_{1,0})
=[u,G_{1,0}]
=\sum_k(k+q)a_{0,k}G_{1,k}.
\end{equation*}
Comparing this with \eqref{eq:even-nonintegral-normalized},
we obtain $p_{r,k}=0$ for $r\ne0$ and
$p_{0,k}=(k+q)a_{0,k}$.

Finally, the bracket formulas give
\begin{equation*}
\begin{aligned}
{}[u,L_{m,i}]&=m\sum_k(k+q)a_{0,k}L_{m,i+k},\\
[u,G_{m,i}]&=m\sum_k(k+q)a_{0,k}G_{m,i+k}.
\end{aligned}
\end{equation*}
These expressions agree with
\eqref{eq:even-nonintegral-normalized} on both bases.
Hence $F=\operatorname{ad}u$, and
$\Delta=D_0+\operatorname{ad}u$ is an even superderivation.
\end{proof}

\subsection{\texorpdfstring{The case $q\in\mathbb Z\setminus\{0\}$}
{The case of nonzero integral q}}
\label{subsec:even-integral}

\begin{lemma}\label{lem:even-local-integral}
If $q\in\mathbb Z\setminus\{0\}$, every local even
superderivation of $S(q)$ is an even superderivation.
\end{lemma}

\begin{proof}
Let $\Delta\in\operatorname{LDer}_{\bar0}(S(q))$.
As in the proof of Lemma~\ref{lem:even-local-nonintegral},
choose an even superderivation $D_0$ at $L_{0,0}+G_{0,1}$ and put
$F=\Delta-D_0$.
Then $F(L_{0,0})=F(G_{0,1})=0$, and
Lemma~\ref{lem:even-basis-action} gives
\eqref{eq:even-nonintegral-normalized} for fixed finite
coefficients $p_{r,k}$.

Apply locality for $F$ at $y=L_{0,0}+G_{1,0}$.
By Proposition~\ref{prop:super-derivations}, its superderivation has
the form
\begin{equation*}
D_y=\operatorname{ad}u+\alpha D_1+\beta D_3+\gamma D_4,
\qquad u=\sum_{r,k}a_{r,k}L_{r,k}.
\end{equation*}
Comparison of parity components gives
$D_y(L_{0,0})=0$ and $D_y(G_{1,0})=F(G_{1,0})$.
Since $D_1(L_{0,0})=D_3(L_{0,0})=0$ and
$D_4(L_{0,0})=qL_{0,-q}$, we obtain
\begin{equation*}
0=D_y(L_{0,0})
=-q\sum_{r,k}r a_{r,k}L_{r,k}+q\gamma L_{0,-q}.
\end{equation*}
The inner term has zero $L_{0,-q}$ coefficient, so
$q\gamma=0$ and hence $\gamma=0$.
The remaining coefficients give $a_{r,k}=0$ for $r\ne0$.
Thus $u=\sum_k a_{0,k}L_{0,k}$.

Now $D_1(G_{1,0})=0$ and $D_3(G_{1,0})=G_{1,-q}$.
Therefore,
\begin{equation*}
F(G_{1,0})
=\sum_k(k+q)a_{0,k}G_{1,k}+\beta G_{1,-q}.
\end{equation*}
Comparing this with \eqref{eq:even-nonintegral-normalized}
gives
\begin{equation*}
p_{r,k}=0\quad(r\ne0),\qquad
p_{0,k}=(k+q)a_{0,k}+\beta\delta_{k,-q}.
\end{equation*}
Put $H=\operatorname{ad}u+\beta D_3$.
The bracket formulas and the definition of $D_3$ yield
\begin{equation*}
H(X_{m,i})
=m\sum_k\bigl((k+q)a_{0,k}+\beta\delta_{k,-q}\bigr)X_{m,i+k},
\qquad X\in\{L,G\}.
\end{equation*}
This agrees with \eqref{eq:even-nonintegral-normalized}
on both bases.
Hence $F=H$, and
$\Delta=D_0+\operatorname{ad}u+\beta D_3$ is an even
superderivation.
\end{proof}

\subsection{\texorpdfstring{The case $q=0$}{The case q = 0}}
\label{subsec:even-zero}

We first establish a result for $\mathcal B(0)$, which will
also be used in the odd case.

\begin{lemma}\label{lem:b0-local}
Every local derivation of $\mathcal B(0)$ is a derivation.
\end{lemma}

\begin{proof}
Put $P(L_{m,i})=\delta_{m,0}\delta_{i,0}L_{0,0}$.
By \cite[Proposition~5]{KaygorodovKhudoyberdiyevShermatova},
every derivation of $\mathcal B(0)$ has the form
$D=\operatorname{ad}u+\alpha D_1+\beta D_2+\rho P$.

Let $A$ be a local derivation.
Locality at $L_{0,0}$ gives $A(L_{0,0})=\rho L_{0,0}$
for some $\rho\in\mathbb C$.
Then $U=A-\rho P$ is local and $U(L_{0,0})=0$.
For $(m,i)\ne(0,0)$, the coefficient of $L_{0,0}$ in
every $D(L_{m,i})$ is zero: the only possible inner
contribution is $[L_{-m,-i},L_{m,i}]=0$.
Thus $U(f)$ has zero $L_{0,0}$ coefficient for every $f$.
If $U(f)=D_f(f)$, comparison of this coefficient gives
$\rho_f f_{0,0}=0$.
Hence the $P$ term vanishes at $f$, and we may write
\begin{equation*}
U(f)=H_f(f),\qquad
H_f=\operatorname{ad}u_f+\alpha_fD_1+\beta_fD_2.
\end{equation*}

For this proof, identify $L_{m,i}$ with $z^mt^i$ in
$R=\mathbb C[z^{\pm1},t^{\pm1}]$.
Write $E_z=z\partial_z$ and $E_t=t\partial_t$.
Then $D_1=E_t$, $D_2=E_z$, and
\begin{equation*}
[u,f]=(E_tu)E_zf-(E_zu)E_tf.
\end{equation*}
Consequently, if $h_z(p)=h_t(p)=0$ at
$p\in(\mathbb C^*)^2$, then $H(h)(p)=0$ for every
$H=\operatorname{ad}u+\alpha D_1+\beta D_2$.

Fix $f\in R$ and $p=(a,b)\in(\mathbb C^*)^2$, and put
\begin{equation*}
h=f-f(p)-f_z(p)(z-a)-f_t(p)(t-b).
\end{equation*}
Since $h_z(p)=h_t(p)=0$, locality gives $U(h)(p)=0$.
Using linearity and $U(1)=0$, we obtain
$U(f)(p)=f_z(p)U(z)(p)+f_t(p)U(t)(p)$.
As Laurent polynomials are determined by their values on
$(\mathbb C^*)^2$, it follows that
\begin{equation*}
U(f)=U(z)\partial_zf+U(t)\partial_tf.
\end{equation*}

Expand the fixed Laurent polynomials as
$z^{-1}U(z)=\sum_{r,k}p_{r,k}z^rt^k$ and
$t^{-1}U(t)=\sum_{r,k}s_{r,k}z^rt^k$.
Both coefficient families have finite support, and
substitution of $f=z^mt^i$ gives
\begin{equation*}
U(L_{m,i})
=\sum_{r,k}(mp_{r,k}+is_{r,k})L_{m+r,i+k}.
\end{equation*}
For $(r,k)\ne(0,0)$, locality at $L_{r,k}$ shows that
the coefficient of $L_{2r,2k}$ is zero.
Indeed, the Euler terms contribute only to $L_{r,k}$,
and the only possible inner contribution is
$[L_{r,k},L_{r,k}]=0$.
Therefore $rp_{r,k}+ks_{r,k}=0$.

Define the finite sum
\begin{equation*}
w=\sum_{\substack{r,k\in\mathbb Z\\k\ne0}}
       \frac{p_{r,k}}{k}L_{r,k}
  -\sum_{r\ne0}\frac{s_{r,0}}{r}L_{r,0}.
\end{equation*}
Using $rp_{r,k}+ks_{r,k}=0$ and
$[L_{r,k},L_{m,i}]=(mk-ri)L_{m+r,i+k}$, we obtain
\begin{equation*}
[w,L_{m,i}]
=\sum_{(r,k)\ne(0,0)}
(mp_{r,k}+is_{r,k})L_{m+r,i+k}.
\end{equation*}
Hence $U=\operatorname{ad}w+s_{0,0}D_1+p_{0,0}D_2$.
Thus $U$, and therefore $A=U+\rho P$, is a derivation.
\end{proof}

\begin{lemma}\label{lem:even-local-zero}
Every local even superderivation of $S(0)$ is an even
superderivation.
\end{lemma}

\begin{proof}
Let $\Delta\in\operatorname{LDer}_{\bar0}(S(0))$.
Its restriction $A$ to $\mathcal B(0)$ is a local derivation,
since every even superderivation restricts to a derivation
of $\mathcal B(0)$.
By Lemma~\ref{lem:b0-local}, $A$ is a derivation.
Write
$A=\operatorname{ad}u+\alpha D_1+\beta D_2+\rho P$.
The even superderivation
$D_0=\operatorname{ad}u+\alpha D_1+\beta D_2+\rho E_L$
extends $A$.
Thus $F=\Delta-D_0$ is local even and
$F(L_{m,i})=0$ for all $m,i$.

Fix $(m,i)\ne(0,0)$.
By locality at $x=L_{m,i}+G_{m,i}$ and
Proposition~\ref{prop:super-derivations}, there is an even
superderivation
\begin{equation*}
D_x=H_x+\rho_xE_L+\sigma_xE_G+\tau_xT,
\qquad
H_x=\operatorname{ad}u_x+\alpha_xD_1+\beta_xD_2,
\end{equation*}
such that $D_x(x)=F(x)$.
The maps $E_L$ and $E_G$ vanish on both input vectors,
and comparison of parity components gives
$H_x(L_{m,i})=0$.
At $q=0$, the inner and Euler terms have identical
coefficients on the corresponding $L$ and $G$ bases.
Hence $H_x(G_{m,i})=0$, and
$F(G_{m,i})=\tau_xG_{m,i}$.
We may therefore write
$F(G_{m,i})=\lambda_{m,i}G_{m,i}$ for every nonzero pair
$(m,i)$.

For distinct nonzero pairs $(m,i)$ and $(n,j)$, apply
locality to
$x=L_{m,i}+L_{n,j}+G_{m,i}+G_{n,j}$.
The same comparison gives
\begin{equation*}
\lambda_{m,i}G_{m,i}+\lambda_{n,j}G_{n,j}
=\tau_x(G_{m,i}+G_{n,j}).
\end{equation*}
Thus all $\lambda_{m,i}$ have a common value $\lambda$.
Locality at $G_{0,0}$ also gives
$F(G_{0,0})=\mu G_{0,0}$ for some $\mu\in\mathbb C$.
Consequently,
$F=\lambda T+(\mu-\lambda)E_G$.
Since both $T$ and $E_G$ are even superderivations,
so is $\Delta=D_0+F$.
\end{proof}

\section{Local odd superderivations}
\label{sec:odd-local-super}

We study local odd superderivations of $S(q)$, treating
the three parameter cases separately.
We first give a normalization for $q\ne0$.

\begin{lemma}\label{lem:odd-normalization}
Let $q\ne0$ and
$\Delta\in\operatorname{LDer}_{\bar1}(S(q))$.
There is $v_0\in S(q)_{\bar1}$ such that
$F=\Delta-\operatorname{ad}v_0$ satisfies
\begin{equation}\label{eq:odd-normalized}
\begin{aligned}
F(G_{0,0})&=0,\qquad F(L_{m,i})=0,\\
F(G_{m,i})&\in
\operatorname{span}_{\mathbb C}\{L_{m,j}:j\in\mathbb Z\}
\end{aligned}
\end{equation}
for all $m,i\in\mathbb Z$.
\end{lemma}

\begin{proof}
Write $\Delta(G_{0,0})=\sum_{r,k}h_{r,k}L_{r,k}$ and put
$v_0=(2q)^{-1}\sum_{r,k}h_{r,k}G_{r,k}$.
Since $[v_0,G_{0,0}]=\Delta(G_{0,0})$, the local odd map
$F=\Delta-\operatorname{ad}v_0$ satisfies $F(G_{0,0})=0$.

By Proposition~\ref{prop:super-derivations}, every odd
superderivation has the form
$D_x=\operatorname{ad}v_x+\delta_xD_6$, where the $D_6$
term is omitted unless $q\in2\mathbb Z\setminus\{0\}$.
Apply locality for $F$ to $x=L_{m,i}+G_{0,0}$ and write
$v_x=\sum_{r,k}a_{r,k}G_{r,k}$.
Since $D_6(G_{0,0})=0$ whenever this map is present,
comparison of parity components gives
\begin{equation*}
0=D_x(G_{0,0})=2q\sum_{r,k}a_{r,k}L_{r,k}.
\end{equation*}
Thus $v_x=0$.
Also, $D_6$ vanishes on the even part, so
$F(L_{m,i})=D_x(L_{m,i})=0$.

Next apply locality to $y=L_{0,0}+G_{m,i}$, with
$v_y=\sum_{r,k}b_{r,k}G_{r,k}$.
Since $F(L_{0,0})=0$ and $D_6(L_{0,0})=0$, we obtain
\begin{equation*}
0=D_y(L_{0,0})
=-q\sum_{r,k}r b_{r,k}G_{r,k}.
\end{equation*}
Hence $b_{r,k}=0$ for $r\ne0$, and
\begin{equation*}
F(G_{m,i})
=2q\sum_k b_{0,k}L_{m,i+k}+\delta_yD_6(G_{m,i}),
\end{equation*}
with the last term omitted when $D_6$ is absent.
The last term can be nonzero only when $m=0$, and its
value then has first index zero.
Thus every term of $F(G_{m,i})$ has first index $m$,
which proves \eqref{eq:odd-normalized}.
\end{proof}

\subsection{\texorpdfstring{The case $q\notin\mathbb Z$}
{The case q not in Z}}
\label{subsec:odd-nonintegral}

\begin{lemma}\label{lem:odd-local-nonintegral}
If $q\notin\mathbb Z$, every local odd superderivation
of $S(q)$ is an inner superderivation.
\end{lemma}

\begin{proof}
Let $\Delta\in\operatorname{LDer}_{\bar1}(S(q))$ and
choose $v_0$ and $F=\Delta-\operatorname{ad}v_0$ as in
Lemma~\ref{lem:odd-normalization}.
Fix $m,i$ and apply locality for $F$ to
$x=L_{1,0}+G_{m,i}$.
Here every odd superderivation is inner, so write
$D_x=\operatorname{ad}v$, where
$v=\sum_{r,k}a_{r,k}G_{r,k}$.
Comparison of parity components gives
\begin{equation*}
F(G_{m,i})=D_x(G_{m,i})
=2q\sum_{r,k}a_{r,k}L_{m+r,i+k},
\qquad D_x(L_{1,0})=0.
\end{equation*}
By \eqref{eq:odd-normalized}, the left-hand side has
first index $m$.
Therefore $a_{r,k}=0$ for $r\ne0$, and
$v=\sum_k a_{0,k}G_{0,k}$.
Using the second equality, we obtain
\begin{equation*}
0=[v,L_{1,0}]
=\sum_k(k+q/2)a_{0,k}G_{1,k}.
\end{equation*}
Since $q/2\notin\mathbb Z$, every factor $k+q/2$ is
nonzero.
Thus $v=0$ and $F(G_{m,i})=0$.
Together with \eqref{eq:odd-normalized}, this gives $F=0$.
Hence $\Delta=\operatorname{ad}v_0$.
\end{proof}

The same proof applies when $q$ is an odd integer, since
$q/2\notin\mathbb Z$ and every odd superderivation is again inner.

\subsection{\texorpdfstring{The case $q\in\mathbb Z\setminus\{0\}$}
{The case of nonzero integral q}}
\label{subsec:odd-integral}

\begin{lemma}\label{lem:odd-local-integral}
If $q\in\mathbb Z\setminus\{0\}$, every local odd
superderivation of $S(q)$ is an odd superderivation.
\end{lemma}

\begin{proof}
For odd $q$, the proof of
Lemma~\ref{lem:odd-local-nonintegral} applies.
Assume that $q$ is even and nonzero.
Let $\Delta\in\operatorname{LDer}_{\bar1}(S(q))$ and
choose $v_0$ and $F=\Delta-\operatorname{ad}v_0$ as in
Lemma~\ref{lem:odd-normalization}.

Whenever we apply locality to $x=L_{0,0}+g$ with
$g\in S(q)_{\bar1}$, its superderivation satisfies
$D_x(L_{0,0})=F(L_{0,0})=0$.
The calculation at $L_{0,0}$ in the proof of
Lemma~\ref{lem:odd-normalization} therefore gives
\begin{equation*}
D_x=\operatorname{ad}\Bigl(\sum_k a_kG_{0,k}\Bigr)
    +\delta_xD_6.
\end{equation*}
Here the coefficients may depend on $x$.

First let $m\ne0$ and apply locality to
$x=L_{0,0}+G_{m,i}+G_{0,0}$.
Since $F(G_{0,0})=0$ and $D_6$ vanishes on both odd
basis vectors, comparison of parity components gives
\begin{equation*}
F(G_{m,i})
=2q\sum_k a_k\bigl(L_{m,i+k}+L_{0,k}\bigr).
\end{equation*}
By \eqref{eq:odd-normalized}, the left-hand side has
first index $m$.
Comparison of the terms with first index zero gives
$a_k=0$ for every $k$.
Hence $F(G_{m,i})=0$ whenever $m\ne0$.
In particular, $F(G_{1,0})=0$.

Now apply locality to
$y=L_{0,0}+G_{0,i}+G_{1,0}$.
Writing the inner element of its superderivation as
$\sum_k b_kG_{0,k}$, we obtain
\begin{equation*}
F(G_{0,i})
=2q\sum_k b_k\bigl(L_{0,i+k}+L_{1,k}\bigr)
 +\delta_y\delta_{i,-3q/2}L_{0,-q}.
\end{equation*}
The left-hand side has first index zero, so comparison
of the terms with first index one gives $b_k=0$ for all $k$.
Therefore $F(G_{0,i})=0$ for $i\ne-3q/2$, while
$F(G_{0,-3q/2})=\delta L_{0,-q}$ for some
$\delta\in\mathbb C$.
Together with \eqref{eq:odd-normalized}, these equalities
give $F=\delta D_6$.
Thus $\Delta=\operatorname{ad}v_0+\delta D_6$ is an odd
superderivation.
\end{proof}

\subsection{\texorpdfstring{The case $q=0$}{The case q = 0}}
\label{subsec:odd-zero}

\begin{lemma}\label{lem:odd-local-zero}
Every local odd superderivation of $S(0)$ is an odd
superderivation.
\end{lemma}

\begin{proof}
Let $\Delta\in\operatorname{LDer}_{\bar1}(S(0))$.
Define the linear bijection
$J:\mathcal B(0)\to S(0)_{\bar1}$ by $J(L_{m,i})=G_{m,i}$,
and put $A(f)=J^{-1}(\Delta(f))$ for $f\in\mathcal B(0)$.

Fix $f\in\mathcal B(0)$.
Let $v_f$ be the inner element of an odd superderivation at $f$,
and put $u_f=J^{-1}(v_f)$.
The bracket formulas and Proposition~\ref{prop:super-derivations}
give
\begin{equation*}
A(f)=[u_f,f]+\alpha_fD_1(f)+\beta_fD_2(f)+\kappa_fP(f),
\end{equation*}
where $P$ is the map used in Lemma~\ref{lem:b0-local}.
The right-hand side is the value at $f$ of a derivation
of $\mathcal B(0)$.
Thus $A$ is a local derivation and, by that lemma,
a derivation.

Write $A=\operatorname{ad}u+\alpha D_1+\beta D_2+\kappa P$,
and define
$D_0=\operatorname{ad}(J(u))+\alpha O_1+\beta O_2+\kappa O_5$.
Then $D_0$ is an odd superderivation agreeing with
$\Delta$ on the even part.
Hence $F=\Delta-D_0$ is local odd and $F(L_{m,i})=0$
for all $m,i$.

At $q=0$, the bracket of any two odd elements is zero.
Consequently, every odd superderivation acts on the odd part as
$\eta\Theta+\tau\Theta_0$.
There are therefore scalars $\lambda_{m,i}$ and $\mu$ such that
\begin{equation*}
\begin{aligned}
F(G_{m,i})&=\lambda_{m,i}L_{m,i},
\qquad (m,i)\ne(0,0),\\
F(G_{0,0})&=\mu L_{0,0}.
\end{aligned}
\end{equation*}
For distinct nonzero pairs $(m,i)$ and $(n,j)$,
locality at $x=G_{m,i}+G_{n,j}$ gives
\begin{equation*}
\lambda_{m,i}L_{m,i}+\lambda_{n,j}L_{n,j}
=\eta_x(L_{m,i}+L_{n,j}).
\end{equation*}
Comparison of coefficients yields
$\lambda_{m,i}=\lambda_{n,j}$.
Denote their common value by $\lambda$.
It follows that $F=\lambda\Theta+(\mu-\lambda)\Theta_0$.
Hence $\Delta=D_0+\lambda\Theta+(\mu-\lambda)\Theta_0$
is an odd superderivation.
\end{proof}

\section{Main result}\label{sec:main-result}

\begin{theorem}\label{thm:local-superderivations}
Every even or odd local superderivation of $S(q)$ is
a superderivation for all $q\in\mathbb C$; that is,
\begin{equation*}
\operatorname{LDer}_{\theta}(S(q))
=\operatorname{Der}_{\theta}(S(q)),
\qquad \theta\in\mathbb Z_2.
\end{equation*}
\end{theorem}

\begin{proof}
Every superderivation is local.
For even maps, the reverse inclusion follows from
Lemmas~\ref{lem:even-local-nonintegral},
\ref{lem:even-local-integral} and \ref{lem:even-local-zero}.
For odd maps, it follows from
Lemmas~\ref{lem:odd-local-nonintegral},
\ref{lem:odd-local-integral} and \ref{lem:odd-local-zero}.
\end{proof}

\end{document}